\documentclass[9pt,twocolumn,twoside]{pnas-new}

\templatetype{pnasresearcharticle} 

\usepackage{amsmath}
\graphicspath{{figures/}}
\usepackage{makecell}

\usepackage{listings} 

\title{Framework for global stability analysis of dynamical systems}

\author[a,1]{George Datseris}
\author[b]{Kalel Luiz Rossi}
\author[c]{Alexandre Wagemakers}

\affil[a]{Department of Mathematics and Statistics, University of Exeter, Exeter, United Kingdom}
\affil[b]{Theoretical Physics/Complex Systems, ICBM, Carl von Ossietzky University Oldenburg, Oldenburg, Lower Saxony, Germany}
\affil[c]{Nonlinear Dynamics, Chaos and Complex Systems Group, Departamento de F\'isica, Universidad Rey Juan Carlos, Madrid, Spain}

\leadauthor{Datseris}

\significancestatement{Dynamical systems are ubiquitous to model natural phenomena ranging from power grids, climate, ecosystems, and more. Global stability analysis is the study of their attracting states over the full state space. It goes beyond (linear) bifurcation analysis and provides insight on the resilience of the states, i.e., how close they are to ``tip’’ to a different, perhaps undesirable stable state such as population death in ecosystems or circulation shutdown in climate. In this work we develop an accurate and flexible framework for global stability analysis of arbitrary dynamical systems, that we believe will accelerate multistability and tipping point research in many fields. We also provide a software implementation that is fast and easy to use in the open source DynamicalSystems.jl library.}

\authorcontributions{Author contributions: G.D. coordinated the project. All authors contributed in compiling the results, creating the code, writing and revising the manuscript.}
\authordeclaration{Authors declare no conflicts of interest.}
\correspondingauthor{\textsuperscript{1}To whom correspondence should be addressed. E-mail: g.datseris@exeter.ac.uk}


\begin{abstract}
Dynamical systems, that are used to model power grids, the brain, and other physical systems, can exhibit coexisting stable states known as attractors. A powerful tool to understand such systems, as well as to better predict when they may ``tip'' from one stable state to the other, is global stability analysis. It involves identifying the initial conditions that converge to each attractor, known as the basins of attraction, measuring the relative volume of these basins in state space, and quantifying how these fractions change as a system parameter evolves. By improving existing approaches, we present a comprehensive framework that allows for global stability analysis on any dynamical system. Notably, our framework enables the analysis to be made efficiently and conveniently over a parameter range. As such, it becomes an essential complement to traditional continuation techniques, that only allow for linear stability analysis. We demonstrate the effectiveness of our approach on a variety of models, including climate, power grids, ecosystems, and more. Our framework is available as simple-to-use open-source code as part of the DynamicalSystems.jl library.
\end{abstract}

\dates{This manuscript was compiled on \today}
\doi{\url{www.pnas.org/cgi/doi/10.1073/pnas.XXXXXXXXXX}}

\begin{document}

\maketitle
\thispagestyle{firststyle}
\ifthenelse{\boolean{shortarticle}}{\ifthenelse{\boolean{singlecolumn}}{\abscontentformatted}{\abscontent}}{}

\dropcap{M}ultistable dynamical systems exhibit two or more co-existing stable states, formally called \emph{attractors}. Multistability is ubiquitous in nature and in mathematical models \cite{feudel2018multistability}, with examples ranging from power grids \cite{hellmann2020network, kim2018multistability, halekotte2021transient}, the climate \cite{marotzke2016instability, lenton2013environmental}, ecosystems like the Amazon rain forest \cite{hirota2011global, dakos2019ecosystem}, the brain and neuronal circuits therein \cite{schwartz2012multistability, kelso2012multistability, kleinschmidt2012variability}, or metabolic systems \cite{zhu2022synthetic, khazaei2022metabolic, geiss2022multistability}. Some attractors of these systems can be desirable, such as synchronized oscillations in power grids, crucial for their proper functioning~\cite{padiyar1999power}. But they can also be undesirable, as for example the extinction of a certain species in ecological models or the collapse of circulation in climate models~\cite{Lohmann2021-wh}.
In a multistable system, the attractor at which the system ends up depends on the initial conditions, but perturbations of the state may enforce switching between attractors, a phenomenon called ``tipping'' ~\cite{Ashwin2012, feudel2018multistability, dakos2019ecosystem}. Alterations in the parameters of a dynamical system can trigger tipping. Hence, it becomes important to evaluate how ``resilient'' attractors are to perturbations, either to parameters or to the system's variables. This is a crucial problem of practical importance in several areas of research \cite{dakos2019ecosystem, halekotte2020minimal}.

A traditional solution to this problem is the \emph{continuation-based bifurcation analysis} (CBA). It identifies fixed points, limit cycles (under some requirements), and describes their linear stability dependence on a system parameter. One major downside is that, by definition, it cannot be applied to chaotic attractors~\cite{DatserisBook}. In the majority of cases, this analysis must be done numerically via one of several software, e.g., AUTO~\cite{doedel1981auto}, MATCONT~\cite{dhooge2003matcont}, CoCo~\cite{Dankowicz2013}, or Bifurcationkit.jl~\cite{veltz:hal-02902346}. The information provided by these frameworks is useful, but incomplete: rigorously speaking, linear stability only conveys information about the system's response to infinitesimally local perturbations. It cannot yield insight on the response to finite perturbations in the state space, which are predominant in practice.

For such responses, it is necessary to study the \emph{global stability}~\cite{Menck2013} of the system's attractors, which involves the nonlinear dynamics over the full state space of the system\footnote{The term ``global stability'' is also used when a dynamical system has a single global attractor, which is different from our use of the term here.}. A proxy for global stability of an attractor is the portion of all possible initial conditions ending up at this attractor, i.e., the \emph{fraction} of the state space that is in the \emph{basin} of said attractor. When the state space is infinite, the concept of the state space fraction becomes a pragmatic one: we need to define a finite-volume box of \emph{physically plausible} initial conditions for the system under study, and we are concerned about the fractions of these plausible initial conditions. In this analysis, attractors with larger basins fractions are globally more stable because stronger perturbations are typically needed to switch the system to another attractor~\cite{Menck2013}. Frequently, this measure is also a much better indicator of the loss of stability as a system parameter is varied, when compared to the linear analysis of the system (see e.g., \cite{Menck2013} or \cite[Chap. 12]{DatserisBook}).

Analyzing global stability as a function of a parameter demands extensive effort from researchers, as it requires the creation of algorithms that can find system attractors and their global stability, ``continue'' them across a parameter, and also perform the expensive numerical simulations required for such algorithms. In the literature, the only framework so far that can aid this analysis is the \emph{featurizing and grouping} approach, proposed first by Gelbrecht et al~\cite{Gelbrecht2020} as MCBB (Monte Carlo Basin Bifurcation Analysis) and then later very similarly by Stender et al~\cite{Stender2021} as bSTAB (basin stability). The method integrates randomly sampled initial conditions (ICs) of a dynamical system for a preset time span. The trajectories of these ICs are then mapped onto \emph{features}, numbers describing the trajectories, such as the mean or standard deviation of some of the system variables. All the feature vectors are clustered using the DBSCAN algorithm~\cite{Ester1996} so that ideally each cluster corresponds to an attractor of the system. The fractions of ICs in each cluster approximate the basin fractions and hence the global stability of each attractor. More details on the method in the Materials and Methods (MatMeth) \S A.

This method works well in a variety of circumstances and can also be applied across a parameter range. However, it comes with significant downsides. One is that it is not clear a-priori which features should be chosen to correctly separate the attractors into clusters, requiring a lot of trial and error. Another downside is that the method cannot guarantee that the clusters of features really correspond to unique attractors, and that it is not mixing two or more attractors together.
An alternative method for finding attractors and their basins of attraction is the \emph{recurrence-based} approach recently proposed in Ref.~\cite{DatserisWagemakers2022}. The method locates attractors by finding recurrences in the system's state space, assuming the Poincar\'e recurrence theorem holds for the system attractors. The input to this method is a state space box, and its tessellation, defining a grid to search for recurrences. Hence, the method will only find attractors within the given box, although the box can initially be arbitrarily large.
We describe the method in more detail in MatMeth\S\textbf{B} and provide a comparison between the two techniques in \S\ref{sec:discussion}\ref{sec:compmeth}. The main advantage of the recurrences method is that it locates the actual system attractors and only requires as an input a state space box that may contain the attractors. So far however it has not been clear how to ``continue'' attractors across a parameter range with this method.

In this work, in \S\ref{sec:results}\ref{sec:recurrences_continuation}, we present a novel global stability analysis and \emph{continuation} algorithm that utilizes the recurrences-based method for finding attractors~\cite{DatserisWagemakers2022}. In \S\ref{sec:results}\ref{sec:applications} we apply it to exemplary models of climate, ecosystem dynamics, and more. As detailed in \S\ref{sec:discussion}, this novel continuation algorithm is the most accurate in finding the actual attractors of a dynamical system, the most transparent in matching attractors across parameter values, and requires the least amount of guesswork from the researcher.
We believe that this novel continuation of global stability, much like global stability analysis itself \cite{Menck2013}, is a crucial new tool for the analysis of dynamical systems. In some cases it supersedes CBA, and in others it can be complemented by CBA, as we discuss in \S\ref{sec:discussion}\ref{sec:compare_linear}.

This continuation method is part of a novel automated framework that performs global stability analysis and continuation, which we present in \S\ref{sec:results}\ref{sec:framework}. Our framework significantly advances existing methodology, including the featurizing methods, thereby including all upsides of current literature while addressing most downsides (MatMeth \S C and D describe the improvements in detail). Its design is based on modular components that can be configured or extended independently. This allows researchers to simply compose the methodology that is best suited to their problem, and then let an automated algorithm execute the process. The framework is accompanied by an extensively tested, well documented, and highly optimized open source software implementation, part of the DynamicalSystems.jl~\cite{DynamicalSystems.jl} general purpose library for nonlinear dynamics (see MatMeth for code example and documentation).

\section{Results}
\label{sec:results}

\subsection{Novel global stability continuation algorithm}
\label{sec:recurrences_continuation}
\begin{figure*}
\centering
\includegraphics[width=17.8cm]{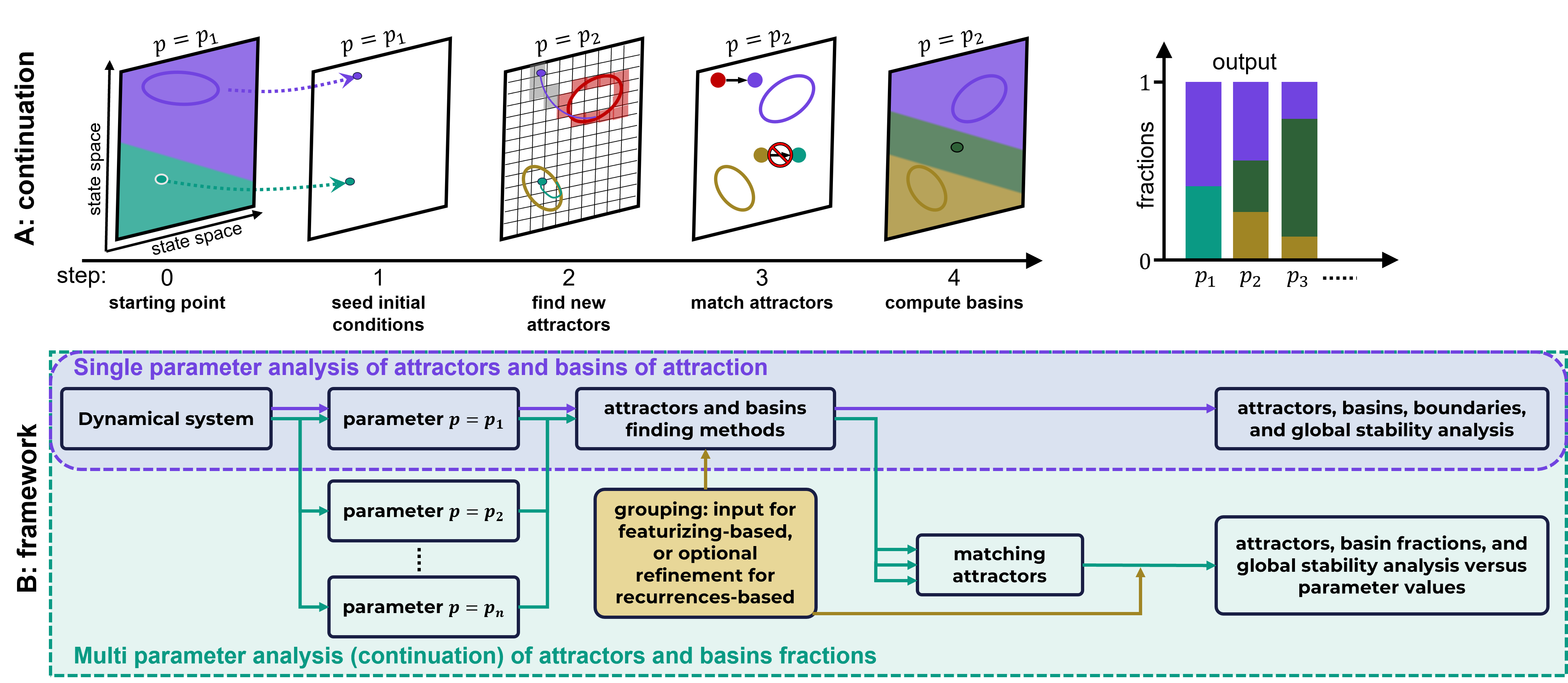}
\caption{\textbf{A}: The recurrences-based seed-and-match algorithm for global stability continuation described in \S\ref{sec:results}\ref{sec:recurrences_continuation}.
\textbf{B}: Schematic illustration of the modular framework for global stability analysis and continuation described in \S\ref{sec:results}\ref{sec:framework}.}
\label{fig:continuation_algorithm}
\end{figure*}

A major contribution of our framework is the novel algorithm for global stability analysis and continuation that we name \emph{recurrences-based attractor find-and-match continuation}, RAFM for short. As illustrated in Fig.~\ref{fig:continuation_algorithm}\textbf{A}, it works as follows.
Step 0: the starting point of the algorithm. Attractors and their basins, or basins fractions, are already known at a parameter $p=p_1$ using the recurrences-based algorithm of Ref.~\cite{DatserisWagemakers2022}.
Step 1: new initial conditions are seeded from the existing attractors. Then, we set the system parameter to $p=p_2$.
Step 2: evolve the seeded initial conditions according to the dynamic rule of the system. The seeds are evolved until they converge to an attractor using the recurrences-based method (the grid reflects the tessellation of the state space that is decided by the user; the finer the grid, the more accurate the results~\cite{DatserisWagemakers2022}).
The main performance bottleneck of the recurrences-based method is finding the attractors. Once found, convergence of other initial conditions is generally much faster~\cite{DatserisWagemakers2022}. To address this in the continuation, we use the observation that, unless a bifurcation is occurring, attractor size, shape, and position, typically depend smoothly on $p$. Hence, the seeded initial conditions at each new parameter will most likely converge the fastest to the new attractors.
Step 3: match attractors in current parameter $p_2$ to those in parameter $p_1$. Matching is arguably the most sophisticated part of the algorithm. It is flexible and attractors can be matched by their ``distance'' in state space, see \S\ref{sec:results}\ref{sec:matching} for more details. In this illustration the ``red'' attractor is matched to the ``purple'' one of Step 0, while the ``yellow'' attractor doesn't map to the previous ``teal'' one, because their state space distance is beyond a pre-defined threshold.
Step 4: with the main bottleneck of the algorithm (finding the attractors) being taken care of, now compute the basins fractions by formally applying the recurrence-based algorithm of Ref.~\cite{DatserisWagemakers2022}. Importantly, the algorithm may still find new attractors (here the ``dark green'' one) that didn't exist before. The end result is the system attractors, and their basins, as functions of the parameter. The attractors and basins are labelled with the positive integers (enumerating the different attractors), and the basins always sum to 1.
Note that because RAFM works on a parameter-by-parameter basis, it can be used to perform continuation across any number of parameters, not just one (we present one here as the simplest conceptual example).

\subsection{Global stability continuation framework}
\label{sec:framework}
To perform global stability analysis, several tasks need to be taken in sequence, see Fig.~\ref{fig:continuation_algorithm}\textbf{B} for an overview. We have abstracted and generalized the tasks to allow researchers different possibilities of how to achieve them.

The first task is the creation of a dynamical system for the global stability analysis. For our framework, this is achieved for free simply by making the implementation part of the DynamicalSystems.jl library~\cite{DynamicalSystems.jl} (see MatMeth \textbf{E}).

The second task is the creation of a mechanism to find attractors and map initial conditions to them. Possibilities for this mechanism are: (1) \emph{featurizing and grouping} initial conditions into attractors (as discussed in the introduction), (2) finding attractors using the \emph{recurrences algorithm} \cite{DatserisWagemakers2022}, or (3) mapping initial conditions to \emph{previously known attractors by proximity}: once the evolution of an initial condition comes close enough to a pre-determined attractor, the initial condition is mapped to that attractor. Note that in (1) or (2) attractors are found via random sampling in the state space, and the probability to find an attractor with basin fraction $f$ after $n$ samples is $1 - (1 - f)^n$.

Mechanism (1) is paired with instructions on how to group features. Currently, the possibilities are: (a) \emph{clustering} features into groups using DBSCAN (as done in MCBB or bSTAB), (b) grouping features by \emph{histograms in the feature space}, so that features that end up in the same histogram bin belong to the same group (novel grouping approach), or (c) mapping features to their \emph{nearest feature} in feature space, from a set of pre-defined features (as done in bSTAB). Having more grouping possibilities than only clustering can be useful, as discussed in \S\ref{sec:results}\ref{sec:matching}.

At this point one can analyze global stability at a given parameter, and also analyze the basin boundaries of attractors for fractal properties. From here, the third task is to ``continue'' the attractors and basins across a parameter range. Our framework has two continuation methods. The first is what has been employed so far by the MCBB or bSTAB algorithms, but with significantly increased accuracy (see MatMeth), and extended to allow any kind of instructions for how to group features. In this approach, trajectories from the dynamical system are generated by sampling random ICs across all parameter values of interest. All these trajectories are mapped to features, that are then grouped, and each group is representing an attractor. The grouped ICs are then re-distributed into the parameter slices they came from, providing the fractions of each group at each parameter value. The alternative is the RAFM algorithm that we described in \S\ref{sec:results}\ref{sec:recurrences_continuation}.

\subsection{Application on exemplary systems}
\label{sec:applications}
\begin{figure}[!t]
    \centering
    \includegraphics[width=8.7cm]{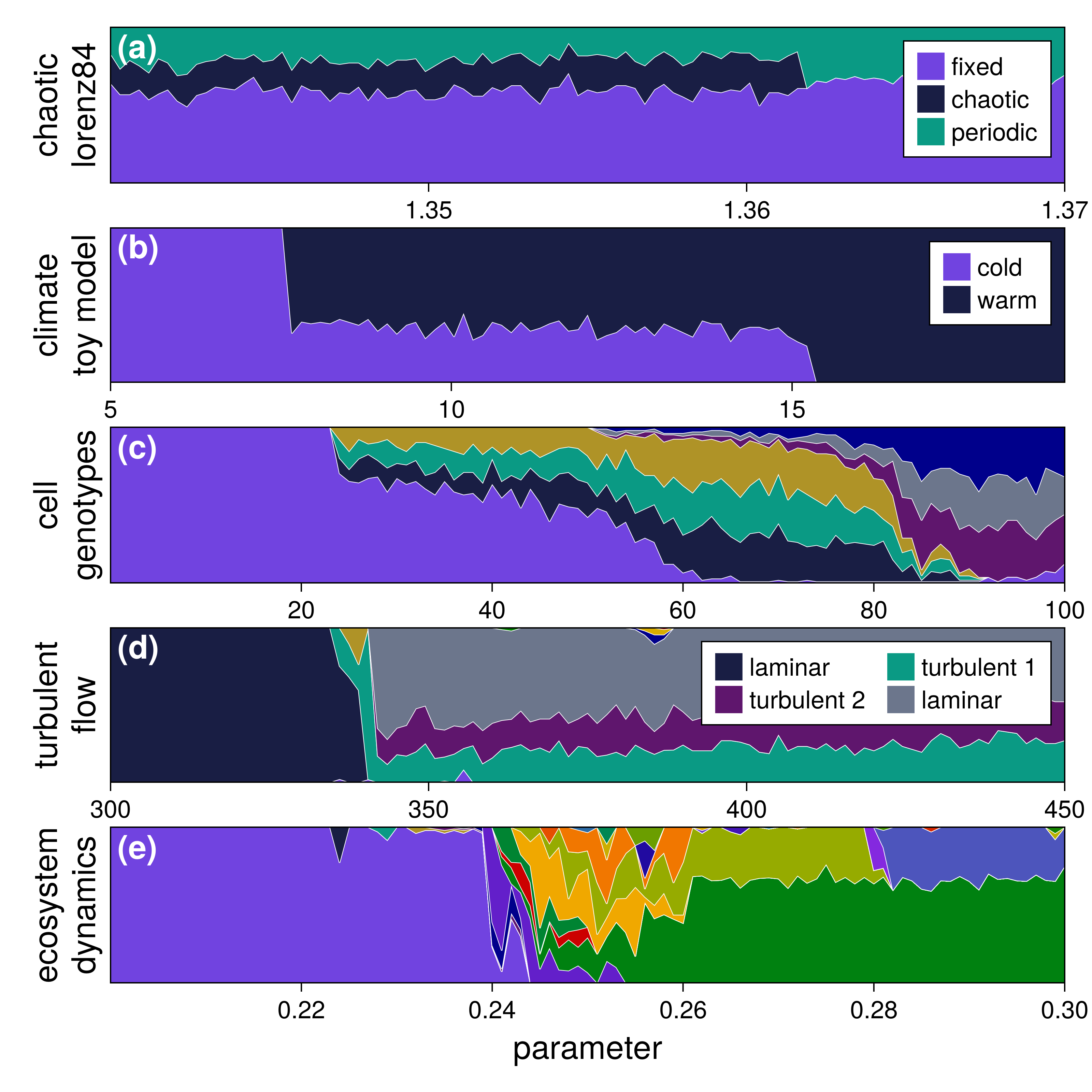}
    \caption{Basins fraction continuation for exemplary dynamical systems using the novel recurrences-based continuation algorithm. The fractions of the basins of attraction are plotted as stacked band plots (hence, summing to 1). Each color corresponds to a unique attractor that is found and continued (but not plotted here). Each simulation scanned 101 parameter values, and in each it sampled randomly 100 initial conditions. The fractions fluctuate strongly versus parameter not due to lack of convergence, but because the basin boundaries are fractal in all systems considered.
    The systems used are:
    a) 3-dimensional paradigmatic chaotic model by Lorenz (Lorenz84~\cite{Lorenz84}) with a co-existence of a fixed point, limit cycle, and chaotic attractor, undergoing a crisis with the chaotic attractor merging into the limit cycle;
    b) 33-dimensional climate toy model~\cite{Gelbrecht2021} featuring bistability of chaotic attractors;
    c) 3-dimensional multistable cell-division model~\cite{huang2006multistability} where each cell type is considered to be a distinct attractor in the gene activity state space;
    d) 9-dimensional model for turbulent shear flow model in which the fluid between two walls experiences sinusoidal body forces~\cite{moehlis2004low};
    e) 8-dimensional ecosystem competition dynamics model~\cite{huisman2001fundamental} featuring extreme multistability (due to the number of attractors, we made no effort to label them further).
    }
    \label{fig:exemplary_fractions}
\end{figure}

In Fig.~\ref{fig:exemplary_fractions} we apply RAFM on some exemplary systems. We stress that we could characterize the different attractors in accurate detail because RAFM finds the actual system attractors, not some incomplete representation of them (i.e., features used in the featurizing-and-grouping approach).

\subsection{Matching and grouping attractors}
\label{sec:matching}
Traditional CBA has a rigid ``matching'' procedure: it always matches the next point found along a ``continuation curve'' to the previous point. This is often correct for infinitesimal perturbations of fixed points, but becomes problematic for global stability analysis, which attempts to find all attractors in a state space box and then continue them. In this case, matching attractors from one parameter to the next becomes a crucial part of the algorithm. For instance, the analysis presented in Fig.~\ref{fig:exemplary_fractions} is only coherent because of the powerful matching procedure implemented in our framework. Without it, the colors would alternate arbitrarily at each parameter value.

In the featurize-and-group algorithm, matching and grouping are the same process. In RAFM, matching operates on a parameter-by-parameter basis. Each time the parameter is incremented, and the new attractors are found, a matching sub-routine is launched. The distance between attractors before and after the parameter change is estimated, with ``distance'' being any symmetric positive-definite function defined on the space of state space sets. By default the Euclidean distance of the attractor centroids is used, but the Hausdorff distance~\cite{Hausdorff1949-nn} is also provided out of the box. After the distance is computed between all new-old attractor pairs, the new attractor labels are matched to the previous attractor labels that have the smallest distance to them, prioritizing pairs with smallest distance. The matching respects uniqueness, so that an attractor from one parameter cannot be matched with more than one attractor from another parameter. Additionally, a distance threshold value can be provided, so that attractors whose distance is larger than this threshold are guaranteed to get assigned different IDs. Note that in principle finding the attractors and matching them are two completely independent processes. If after the continuation process is finished the user decides that the chosen matching procedure was unsuitable, they can launch a ``re-matching'' algorithm with different matching ``distance'' function, without having to re-do any computations for finding the attractors or their fractions (matching only renames attractor labels).

\begin{figure}
    \centering
    \includegraphics[width=8.7cm]{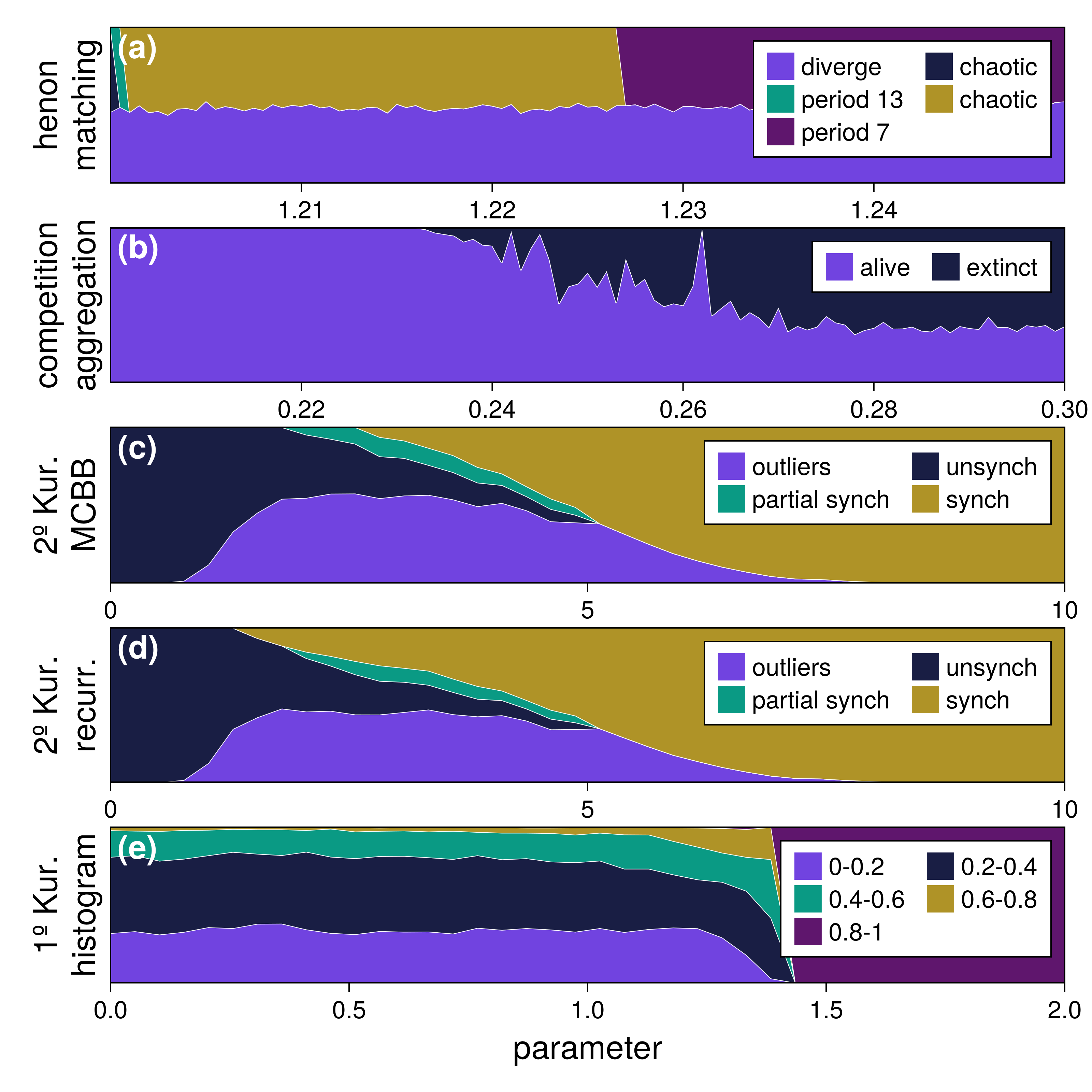}
    \caption{Highlights of the matching or grouping components of the framework. a) Matching attractors of the H\'enon map based on their period. In the chosen parameter range, an attractor is transformed from chaotic, to period 3, 7, and 14. The attractor stays in approximately the same state space location, so whether we consider the centroid distance or the Hausdorff distance, the attractor would be matched to itself in all parameter values due to the very small distance evaluation. Here however we use as distance $f(A,B) = |\log_2(\mathrm{len}(A)) - \log_2(\mathrm{len}(B))|$, with $\mathrm{len}$ measuring the amount of cells (of the state space tessellation) the attractor covers, and threshold $t=0.9\bar{9}$. This effectively means that matched attractors must have periods with ratio less than 2. b) Grouping attractors of Fig.~\ref{fig:exemplary_fractions}(e) so that attractors are grouped into those whose 3rd species has population less than 0.01, or more. (c) A replication of the MCBB~\cite{Gelbrecht2020} results for a 2nd-order Kuramoto oscillator network representing a power grid, using the featurize-and-group continuation implementation from our framework. Features extracted from sampled trajectories are the means of the frequencies. (d) Same system as (c), but using the recurrence continuation and matching attractors by their centroid distance (i.e., as in Fig.~\ref{fig:exemplary_fractions}). The only extra step was to post-process the results so that all attractors with basins fractions less than 4\% are aggregated (as was done in Ref.~\cite{Gelbrecht2020} and in panel (c)). (e) Attractor basin fractions for a network of 1st order Kuramoto oscillators; the attractors here are found and matched using the recurrences continuation, and then grouped via a histogram of their synchronization order parameter $R$~\cite[Chap. 9]{DatserisBook}. Attractors whose order parameter $R$ falls in the same histogram bin are aggregated.
    }
    \label{fig:matching}
\end{figure}

The last thing to highlight in this section is the desirable post-processing of grouping similar enough attractors. This happens automatically if one uses the featurize-and-group continuation method. However, taking as an example Fig.~\ref{fig:exemplary_fractions}(e), the RAFM method finds countless individual attractors. For the researcher, the individual attractors may be useful for careful analysis, but it is sometimes desirable to group similar enough attractors. In our framework it is possible to use exactly the same grouping infrastructure utilized by the featurizing-and-grouping continuation, but now applied to the outcome of RAFM as a post-processing step. In Fig.~\ref{fig:matching} we highlight examples that utilize the powerful matching and/or grouping components offered by our framework.

\section{Discussion}
\label{sec:discussion}

\subsection{Comparison with traditional continuation-based bifurcation analysis}
\label{sec:compare_linear}
\begin{figure*}
\centering
\includegraphics[width=17.8cm]{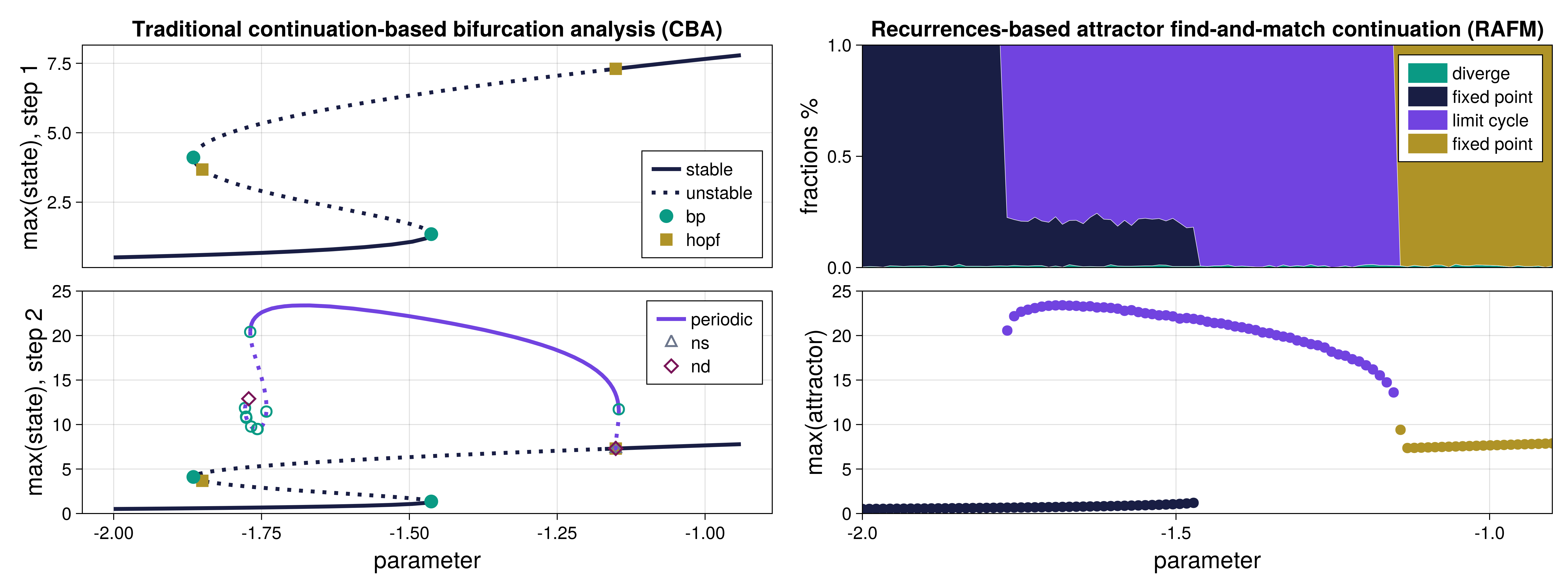}
\caption{Stability analysis of a 3-dimensional neural dynamics model~\cite{Cortes2013}, plotting as information the maximum of the first variable of the dynamical system. For the parameters used, the model undergoes saddle-node and Hopf bifurcations and features bistability of a limit cycle and a fixed point. Left: analysis using the BifurcationKit.jl~\cite{BifurcationKit.jl} Julia package for CBA. The analysis must happen in two steps; first the branch of a fixed point is found and continued, and then, using a different algorithm, the branch of the limit cycle is continued from the Hopf bifurcation. We scatter-plot special points found and labelled by the process. The computation required $\sim$19 seconds on an average laptop. Right: analysis of the same model using our framework. The analysis happens in one fully automated step, after deciding the state space box and other meta-parameters. For the analysis, we purposefully demanded unnecessarily high accuracy, using a 9-th order ODE solver~\cite{Verner2010, DifferentialEquations.jl-2017} with tolerances of 10\textsuperscript{-9}, and requiring 1000 recurrences before claiming convergence in the recurrence algorithm~\cite{DatserisWagemakers2022}. The process integrated in total 101,000 initial conditions, yet required $\sim$16 seconds on the same laptop. Attractor matching utilized a threshold: attractors whose distance in terms of their maximum value of the first variable (i.e., same information plotted in figure) exceeding 3.0 are not matched.}
\label{fig:comparison_with_bifkit}
\end{figure*}

In Table~\ref{tab:comparison} and Fig.~\ref{fig:comparison_with_bifkit} we provide a careful comparison between CBA and RAFM. A direct comparison of the two approaches is difficult, since their operational basis, and main output, are fundamentally different. On one hand, CBA finds and continues the curves of individual fixed points or limit cycles across the joint state-parameter space on the basis of Newton's method. On the other hand, RAFM first finds attractors at all parameter values using the original system equations, and then \emph{matches} appropriately similar attractors across different parameters, giving the illusion of continuing them individually. Additionally, the curves of stable fixed points in the joint parameter space (Fig.~\ref{fig:comparison_with_bifkit}) are only a small part of the information our framework provides. Important provided information are the basin fractions and how they change, which is completely absent in CBA.

Based on this comparison, we argue that the RAFM algorithm and the global stability framework we provide is an essential tool for stability analysis of dynamical systems. We further believe that in some application scenarios RAFM will supersede CBA, especially given the difference in required user expertise, required interventions, and steepness of the learning curve that CBA has over RAFM. In other scenarios, we envision that RAFM can be used as the default analysis method, providing the majority of information, and CBA then becomes a more in-depth analysis of fixed points, limit cycles, and bifurcations, if such analysis is required.

\newcounter{tablecounter} 

\begin{table*}[]
    \centering
    \begin{tabular}{lll}\toprule
    \textbf{\#} & \textbf{Traditional continuation-based bifurcation analysis (CBA)} & \textbf{Recurrences-based attractor find-and-match continuation (RAFM)} \\ \hline
    \rowcolor{blue!20} \refstepcounter{tablecounter}\arabic{tablecounter} & provides curves of unstable fixed points / limit cycles & only finds attracting, not repelling, sets \\
    \rowcolor{blue!20} \refstepcounter{tablecounter}\arabic{tablecounter} & several possibilities for how to continue bifurcation curves & does not explicitly detect bifurcation points \\
    \rowcolor{blue!20} \refstepcounter{tablecounter}\arabic{tablecounter} & does not put limits on state space extent & needs as an input a state space box that may contain attractors \\
    \rowcolor{blue!20} \refstepcounter{tablecounter}\arabic{tablecounter}$^\diamond$ & likely to find fixed points / limit cycles with small or even zero basins & probability to find attractor is proportional to its basins fraction \\
    \rowcolor{blue!20} \refstepcounter{tablecounter}\arabic{tablecounter} & detects and classifies local bifurcation points & does not compute Jacobian eigenvalues at all \\
    \hline
    \rowcolor{green!20} \refstepcounter{tablecounter}\arabic{tablecounter} & finds and continues fixed points and periodic orbits & finds and continues any kind of attractors, including quasiperiodic or chaotic \\
    \rowcolor{green!20} \refstepcounter{tablecounter}\arabic{tablecounter} & requires a computable Jacobian of the dynamic rule & works for any dynamical system including Poincar\'e maps or projections \\
    \rowcolor{green!20} \refstepcounter{tablecounter}\arabic{tablecounter}$^\mathparagraph$ & user must manually search for multistability & different attractors are automatically detected (via random sampling) \\
    \rowcolor{green!20} \refstepcounter{tablecounter}\arabic{tablecounter} & does not compute the basins of attraction or their fractions & computes the fractions and, if computationally feasible, also the full basins \\
    \rowcolor{green!20} \refstepcounter{tablecounter}\arabic{tablecounter} & uses the local, linearized dynamics & uses the full nonlinear dynamics \\
    \rowcolor{green!20} \refstepcounter{tablecounter}\arabic{tablecounter}$^\ddagger$ & limited use in indicating loss of stability & more likely to indicate loss of stability as the basin fraction approaches 0 \\
    \rowcolor{green!20} \refstepcounter{tablecounter}\arabic{tablecounter}$^*$ & parameter change may not affect linear stability of all fixed points & parameter change is more likely to affect global stability of all attractors \\
    \rowcolor{green!20} \refstepcounter{tablecounter}\arabic{tablecounter} & no sophistication on matching fixed points & sophisticated, user-configurable matching of attractors \\
    \rowcolor{green!20} \refstepcounter{tablecounter}\arabic{tablecounter}$^\dagger$ & requires expertise and constant interventions  & conceptually straightforward even in advanced use-cases \\
    \bottomrule\end{tabular}
    \caption{A comparison between CBA and RAFM as tools for analyzing the stability of a dynamical system versus a parameter. The table is colored blue or green for when to prefer CBA or RAFM respectively.}
    \addtabletext{
    $^\diamond$Newton's method transforms the dynamical system into a discrete time system with different basins, making attractors with very small or zero basin sizes have much larger ones instead.
    $^\mathparagraph$In RAFM, attractors that are not being continued from a previously found one, are found via random sampling of initial conditions in the given state space box. The probability to find an attractor is equal to $1 - (1 - f)^n$ with $f$ the basins fraction of the attractor and $n$ the amount of sampled initial conditions.
    $^\ddagger$Changing a parameter often does not meaningfully increase the unstable eigenvalues of the Jacobian matrix, which would indicate loss of stability ~\cite[Chap. 12]{DatserisBook}. On the other hand, basin fractions typically decrease smoothly towards zero as an attractor loses stability~\cite{Menck2013}, although this is not guaranteed to be the case~\cite{Schultz2017}, in which scenario, neither method indicates loss of stability.
    $^*$Change of a parameter may affect the linear stability of a single fixed point, not all, providing flat lines in the bifurcation diagram for the unaffected fixed points. On the contrary, loss of global stability of any attractor affects (typically increases) the global stability of all other attractors.
    $^\dagger$Advanced applications of traditional bifurcation analysis software require several manual interventions during the process, and tuning of several configuration options, many of which do not have an immediately transparent role, requiring an expert user to make several decisions. The simplicity of our approach comes in part because of the brute-force nature of mapping individual initial conditions to attractors to collect the fractions, the intuitive nature of how attractors are matched (which is also user configurable), and the lack of necessity of interventions: after the configuration is decided, the framework runs automatically.}
    \label{tab:comparison}
\end{table*}

\subsection{Comparison between attractor-finding methods}
\label{sec:compmeth}
Our framework provides two radically different methods for finding attractors: the recurrence-based and the featurize-and-group approach. Generally speaking, the recurrence-based method should be preferred when possible, because of its accuracy (finding the actual attractors), and the possibility for follow-up analysis of found attractors. However, the featurize-and-group method should be preferred when the recurrence-based method fails, because e.g., the provided state space tessellation is ill-defined, or because computational demands exceed what is available. In Table~\ref{tab:comparison_attractors} we provide a comprehensive comparison between the two methods.

 \begin{table*}[!t]
     \centering
     \begin{tabular}{p{1.5cm}p{7.5cm}p{7.5cm}}\toprule
     \textbf{Aspect} & \textbf{Recurrences method} & \textbf{Featurizing method} \\ \hline

     Accuracy & Highly accurate: finds actual attractors using their unique property of their state space location. &
     Less accurate, as trajectories are transformed into features, and attractors correspond to groups of features. The correspondence is not guaranteed to be unique or reversible. \\

     Info stored & Stores samples of points on the found attractors & Stores a user-provided function of the group of features (by default: the centroid of the group). \\

     Speed & Very fast for low dimensional systems and for systems whose attractors are fixed points or periodic orbits. Becomes slow for attractors with long recurrence times (high dimensional chaotic attractors or very fine state space tessellations). & Performance is independent of system attractors. It linearly scales with the amount of initial conditions, the transient integration time, and the total integration time. Parallelizable. In addition is the cost of the grouping process, which is huge for clustering but trivial for histograms or nearest-feature. See also the benchmark comparison in the Supplementary Information \S 3 \\

     Memory & Memory allocation scales as $(1/\varepsilon)^\Delta$, with $\varepsilon$ the state space tessellation size and $\Delta$ the capacity dimension of the attractor, which is often much lower than the state space dimension~\cite{DatserisBook}. & Memory allocation of the trajectories scales linearly with integration time and sampling rate. Additionally, specifically for clustering used as the grouping mechanism, the total memory allocated is proportional to the square of [initial conditions $\times$ parameter values] which, if one attempts to obtain accurate results, is often beyond the available memory on a typical computer.\\

     Necessary input (guesswork) & A state space box that may contain the attractors, and a state space tesselation that is fine enough to differentiate the location of attractors. & A state space box that may contain the attractors; a function mapping attractors to features, such that different attractors produce different features; how much transient time to discard from time evolution; how much time to evolve and record the trajectory for, after transient. \\

     Meta-parameters & The parameters of the finite state machine of the recurrences algorithm discussed in Ref.~\cite{DatserisWagemakers2022}, and the time stepping $\Delta t$. All are crucial, but all are conceptually straightforward. & The integration time and sampling rate and all parameters of the grouping procedure. Integration parameters are straightforward, but optimal parameters related to the grouping are much harder to guess.\\

     Troubleshooting & Easy to troubleshoot. At any point the actual trajectories and attractors are accessible and why a failure occurs is typically easy to find out. Matching of attractors happens parameter-by-parameter, hence, individual parameter slices can be isolated and analyzed to identify where matching has failed and why (distances between attractors is also provided information). & Difficult to troubleshoot. It does not find the actual system attractors, so the user must always reason in terms of features. When grouping using clustering (DBSCAN), the grouping process essentially operates as a black box after the features have been computed, making it harder to comprehend failures. Matching of attractors happens at the same time as grouping, making it nearly impossible to understand why an expected matching failed during the continuation process. \\

     Failures & Algorithm may fail if state space tessellation is not fine enough and a grid cell may contain points from different attractors. Chaotic saddles and other sticky sets generate very long transients \cite{lai2009transient} and the algorithm can interpret them as attractors. Additionally, the algorithm is sensitive to the time step $\Delta t$ and the used integrator. For limit cycles it is often the case that a non-adaptive integrator needs to be used. & Sticky sets that are formally not attractors will be interpreted as attractors. Additionally, for ill-defined features, any group of trajectories could be interpreted as an attractor, which is not desirable in the context of this paper. Clustering via DBSCAN may fail unexpectedly, or, finding the optimal radius for the clustering may yield incorrect results. \\
     \bottomrule
     \end{tabular}
     \caption{Comparison the two main methods of finding and continuing attractors.}
     \label{tab:comparison_attractors}
 \end{table*}

\subsection{Future extensions}
The modularity of our framework allows it to be easily extended by researchers for their specific use cases. For instance, the default information about the attractors that is ``tracked'' during the RAFM continuation is their state space location. In the future, we plan to further enhance the framework to also (optionally) track dynamical characteristics of the attractors such as the Lyapunov spectra or fractal dimensions~\cite{DatserisBook} (computing these quantities after-the-fact is trivial as illustrated in our example code of Listing.~\ref{lst:code}). Utilizing this information will allow us to better label the attractors (instead of using the integers as we currently do), better match attractors during the continuation, or even automatically detect bifurcations or changes in the nature of the attractor, such as transitioning from a periodic attractor (no positive Lyapunov exponents) to a chaotic one (at least one positive Lyapunov exponent).

\matmethods{
\label{sec:methods}

\newcounter{matmethcounter}

\subsection*{\refstepcounter{matmethcounter}\Alph{matmethcounter}. Featurizing methods for finding attractors}
\label{mat:featurizing_summary}
We group together two similar methods that have been recently proposed in the literature for finding attractors. One is called Monte Carlo basin bifurcation analysis (MCBB) \cite{Gelbrecht2020}, the other is basin stability analysis (bSTAB) \cite{Stender2021}. Both methods work by identifying attractors as clusters of user-defined-features of trajectories. Their first step is to integrate $N$ randomly chosen initial conditions inside a certain box in state space. The integration is done for some time $T$, after a transient $T_\mathrm{tr}$. Both $T$ needs to be sufficiently large so that the trajectories correspond to the systems' long-term behavior and also $T_\mathrm{tr}$ needs to be large to avoid the transient regime. Each trajectory $\vec{x}(t)$ is then transformed into a vector of $K$ features $\vec{F}$, specified by some featurizer function $\vec{f}$ such that $\vec{f}(\vec{x}(t)) = \vec{F}$ that has to be defined by the user. Each vector of features $\vec{F}$ describes a point in the $K$-dimensional space of features $f_1 \times f_2 \times \cdots \times f_K$. The key idea is that features belonging to the same attractor cluster together in state space, so that each attractor forms a distinct cluster (a cloud of points) in feature space. The final step in the method is to therefore cluster the features. The clustering algorithm chosen for this is the Density-based Spatial Clustering of Applications with Noise (DBSCAN) \cite{Ester1996}. It first classifies two points as neighbors if their distance if smaller than a radius $\epsilon$. Then, it clusters together points with many neighbors (equal or more than a parameter minPts), and leave as outliers points with too few neighbors (less than minPts). The radius $\epsilon$ is a crucial parameter for the algorithm, and often needs fine tuning for proper clustering. The methods by \cite{Gelbrecht2020} and \cite{Stender2021} use two different ways to identify a value for $\epsilon$. Authors in \cite{Gelbrecht2020} use a method that looks at the ordered distance of the $k$ nearest neighbors to each point in the dataset, and finds the first knee (high derivative point) \cite{Ester1996, Hahsler2019}. Authors in \cite{Stender2021} iteratively search for the $\epsilon$ that maximizes a criterion of clustering quality. To calculate this criterion, they evaluate the silhouette of the each cluster, which measures how similar each point is to the cluster it currently belongs to, compared to the other clusters, and ranges from $-1$ (worst matching) to $+1$ (ideal matching). This leads to one silhouette value per feature; the authors then take the minimum value as the representative for the clustering quality for each radius. The chosen radius is thus the value that leads to the highest minimum silhouette. In both methods, the clusters found by DBSCAN are considered then as attractors.

\subsection*{\refstepcounter{matmethcounter}\Alph{matmethcounter}. Recurrences-based method for finding attractors}
\label{mat:recurrences_summary}
The inputs to this method are a dynamical system, a state space box that may contain the attractors (although initially it may be arbitrarily large), and a tessellation of the given box into cells. An initial condition of the system is evolved step-by-step with time step $\Delta t$. At each step, the location of the trajectory in state space is mapped to its cell, and that cell is labelled as visited. If the dynamical system has attractors, and they satisfy the Poincar\'e recurrence theorem~\cite[Chap. 9]{DatserisBook}, the trajectory is guaranteed to revisit cells it has visited before. Once a pre-decided number of recurrences $n_f$ have been accumulated consecutively (i.e., enough previously-visited cells are visited again), the method claims to have found an attractor. It then proceeds to locate the attractor accurately, by collecting a pre-decided number $n_l$ of points on the attractor (Figure 1 of Ref.~\cite{DatserisWagemakers2022} and panel (3) of Fig.~\ref{fig:continuation_algorithm}). A finite state machine formulation keeps track of coexisting attractors, so that each attractor is unique. It also keeps track of divergence to infinity by counting steps $n_d$ outside the box, ensures algorithm termination by setting a total $n_m$ of max amount of $\Delta t$ iterations, and makes convergence faster by utilizing information already encoded in the grid: if the trajectory visits consecutively a relatively small number $n_r$ of cells already labelled as an attractor, convergence is already eagerly decided. I.e., converging to an already found attractor is much faster than finding that attractor for the first time. Hence, $\Delta t, n_f, n_l, n_d, n_m, n_r$ are the meta-parameters of the algorithm and have sensible default values that work in most cases. More information on the method can be found in Ref.~\cite{DatserisWagemakers2022}. Notice that the recurrence method is different at a fundamental level from GAIO\cite{gerlach2020set} and other cell mapping techniques\cite{sun2018cell}. We expand more on this in the Supplementary Information. Also note that the method is not perfect; it may identify two attractors when only one exists due to e.g., choosing too low convergence criteria $n_l, n_f$, or due to a commensurate period of attractor and integrator time step. However, once again the importance of finding the ``actual'' attractors becomes apparent: further analysis, by e.g., plotting the attractors, immediately highlights such a failure and how to deal with it, and we also provide several tips in the documentation of our method in the code implementation~\cite{Attractors.jl}.

\subsection*{\refstepcounter{matmethcounter}\Alph{matmethcounter}. Improvements to the recurrences method}

A large drawback of the recurrences method was that it scaled poorly with the dimension $D$ of the dynamical system. If an $\epsilon$-sized tessellation of the state space is chosen, then memory allocated scaled as $1/\epsilon^D$. We now use sparse arrays to store accessed grid locations. This changes the memory scaling to $1/\epsilon^\Delta$, with $\Delta$ the capacity dimension~\cite{DatserisBook} of the attractor, which is typically much smaller than $D$ (and is only $1$ for limit cycles).

\subsection*{\refstepcounter{matmethcounter}\Alph{matmethcounter}. Improvements to the featurizing methods}
First, we have changed the criterion used for finding the optimal radius in the clustering method. We have found the knee method consistently more unreliable than the iterative search. We have also found that the mean, instead of the minimum, silhouette value as the measure of clustering quality leads to better clustering. For instance, this lead to correct clustering in the Lorenz86 system, whereas the minimum value criterion did not. Furthermore, our method searches for the optimal-radius radius with a bisection method, instead of the linear method used by authors in \cite{Stender2021}. This significantly speeds up the code.
Another simple modification we introduced is to rescale the features in each dimension ($f_1, f_2, \cdots, f_K$) into the same interval, for instance $[0,1]$. We noticed that the clustering method performs poorly if the features span different ranges of values, and this simple modification proved to be a very powerful solution.
Third, we allow the integration of all initial conditions to be done in parallel, using several computer cores, which speeds up the solution.
Lastly, grouping in our framework can also happen based on a histogram in feature space.

\subsection*{\refstepcounter{matmethcounter}\Alph{matmethcounter}. Code implementation}

The code implementation of our framework is part of the DynamicalSystems.jl library~\cite{DynamicalSystems.jl} as the Attractors.jl package~\cite{Attractors.jl}. The code is open source code for the Julia language, has been developed following best practices in scientific code~\cite{goodscientificcode}, is tested extensively, and is accompanied by a high quality documentation. An example code snippet is shown in Listing~\ref{lst:code}.

Besides the quality of the implementation, three more features of the code are noteworthy. First, that it is part of DynamicalSystems.jl instead of an isolated piece of code. This integration makes the simplicity and high-levelness of Listing~\ref{lst:code} possible, and makes the input for the code easy to set up. Moreover, the direct output of the code can be used with the rest of the library to further analyze attractors in terms of e.g., Lyapunov exponents or fractal dimensions. Indeed, in the provided code example Listing~\ref{lst:code} we compute the Lyapunov spectra of all found attractors, across all parameter values, in only two additional lines of code. Second, utilizing the Julia language's multiple dispatch system~\cite{bezanson2017julia}, the code is extendable. It establishes one interface for how to map initial conditions to attractors, and one for how to group features, both of which can be extended, and yet readily be usable by the rest of the library such as the continuation methods. Third, a lot of attention has been put into user experience, by establishing a short learning curve via a minimal user interface, by carefully considering how to provide the output in an intuitive format, as well as providing easy-to-use plotting functions that utilize the code output. More overview and information on the code or its design can be found in its online documentation or source code~\cite{Attractors.jl}.

\subsection*{\refstepcounter{matmethcounter}\Alph{matmethcounter}. Code for this article}
The code we used to create the figures of this article is fully reproducible and available online \cite{codebase}.

} 

\showmatmethods{} 

\begin{lstlisting}[language=Python, label = {lst:code},basicstyle=\footnotesize\ttfamily, caption = {Julia code snippet showcasing the usage of the DynamicalSystems.jl implementation of our framework. The code produces panel (a) of Fig. 2. The main output of the code are two vectors, containing the basins fractions and attractors at each parameter value respectively. The fractions and attractors are formulated as dictionaries, mapping attractor labels (the integers) to basin fractions and sets of points on the attractor, respectively. At its end, the code snippet computes the Lyapunov spectra of all found attractors, by using the first point on each attractor as initial condition for the computation of the Lyapunov spectrum.}]
using DynamicalSystems # our framework implementation
using OrdinaryDiffEq   # high-accuracy ODE solvers

# create Lorenz84 within DynamicalSystems.jl
function lorenz84_rule(u, p, t)
    F, G, a, b = p
    x, y, z = u
    dx = -y^2 -z^2 -a*x + a*F
    dy = x*y - y - b*x*z + G
    dz = b*x*y + x*z - z
    return SVector(dx, dy, dz)
end
u0 = ones(3) # init. state
p0 = [6.886, 1.347, 0.255, 4.0] # init. parameters
# ODE solver:
diffeq = (alg = Vern9(), reltol = 1e-9, abstol = 1e-9)
# Main object of the library: a `DynamicalSystem`
ds = CoupledODEs(lorenz84_rule, u0, p0; diffeq)

# Provide state space box tessellation to search in
xg = yg = zg = range(-3, 3; length = 600)
grid = (xg, yg, zg)
# initialize recurrences-based algorithm
# and choose its metaparameters
mapper = AttractorsViaRecurrences(ds, grid;
    mx_chk_fnd_att = 1000, mx_chk_loc_att = 2000,
    mx_chk_lost = 100, mx_chk_safety = 1e8,
    Dt = 0.05, force_non_adaptive = true,
)

# find and continue attractors across a given
# parameter range for the `pidx`-th parameter
prange = range(1.34, 1.37; length = 101)
pidx = 2 # index of parameter
sampler = statespace_sampler(grid)[1]
rsc = RAFM(mapper)
# main output:
fractions_curves, attractors_info = continuation(
    rsc, prange, pidx, sampler
)

# Estimate Lyapunov spectra for all attractors
# by looping over the parameter range
lyapunovs_curves = map(eachindex(prange)) do index
    set_parameter!(ds, pidx, prange[index])
    attractor_dict = attractors_info[index]
    exponents = Dict(
        k => lyapunovspectrum(ds, 10000; u0 = A[1])
        for (k, A) in attractor_dict
    )
end
\end{lstlisting}

\acknow{
The authors would like to thank Ulrike Feudel, Peter Ashwin, Ulrich Parlitz, and Harry Dankowicz for helpful discussions.
G.D. was supported by the Royal Society International Newton Fellowship.
K.L.R. was supported by the German Academic Exchange Service (DAAD).
A.W. was supported by the Spanish State Research Agency (AEI) and the European Regional Development Fund.
}

\showacknow{} 
\bibliography{REFERENCES}

\end{document}